\documentclass[12pt]{article}%
\usepackage{graphicx}
\usepackage[intlimits]{amsmath}
\usepackage{latexsym}
\usepackage{amsfonts}
\usepackage{amssymb}%
\makeatletter
\newfont{\footsc}{cmcsc10 at 8truept}
\newfont{\footbf}{cmbx10 at 8truept}
\newfont{\footrm}{cmr10 at 10truept}
\newtheorem{theorem}{Theorem}

\newtheorem{claim}[theorem]{Claim}

\newtheorem{corollary}[theorem]{Corollary}

\newtheorem{lemma}[theorem]{Lemma}

\newenvironment{proof}[1][Proof]{\noindent{\textbf {#1}  }}  {\hfill$\Box$\bigskip}

\begin{document}

\title{Degree powers in graphs with forbidden even cycle}
\author{Vladimir Nikiforov\\Department of Mathematical Sciences, University of Memphis, \\Memphis TN 38152, USA, email: \textit{vnikifrv@memphis.edu}}
\maketitle

\begin{abstract}
Let $C_{l}$ denote the cycle of length $l.$ For $p\geq2$ and integer $k\geq1,$
we prove that the function
\[
\phi\left(  k,p,n\right)  =\max\left\{  \sum_{u\in V\left(  G\right)  }%
d^{p}\left(  u\right)  :G\text{ is a graph of order }n\text{ containing no
}C_{2k+2}\right\}
\]
satisfies $\phi\left(  k,p,n\right)  =kn^{p}\left(  1+o\left(  1\right)
\right)  .$

This settles a conjecture of Caro and Yuster.

Our proof is based on a new sufficient condition for long paths, that may be
useful in other applications as well.

\end{abstract}

\section{Introduction}

Our notation and terminology follow \cite{Bol98}; in particular, $C_{l}$
denotes the cycle of length $l$.

For $p\geq2$ and integer $k\geq1,$ Caro and Yuster \cite{CaYu00} studied the
function
\[
\phi\left(  k,p,n\right)  =\max\left\{  \sum_{u\in V\left(  G\right)  }%
d_{G}^{p}\left(  u\right)  :G\text{ is a graph of order }n\text{ without a
}C_{2k+2}\right\}
\]
and conjectured that
\begin{equation}
\phi\left(  k,p,n\right)  =kn^{p}\left(  1+o\left(  1\right)  \right)  .
\label{YC}%
\end{equation}

The graph $K_{k}+\overline{K}_{n-k},$ i.e., the join of $K_{k}$ and
$\overline{K}_{n-k},$ gives $\phi\left(  k,p,n\right)  >k\left(  n-1\right)
^{p},$ so to prove (\ref{YC}) a matching upper bound is necessary. We give
such a bound in Corollary \ref{cor1} below. Our main tool, stated in Lemma
\ref{le1}, is a new sufficient condition for long paths. It also implies the
following spectral bound proved in \cite{Nik08}:\medskip

\emph{Let }$G$\emph{ be a graph of order }$n$\emph{ and }$\mu$\emph{ be the
largest eigenvalue of its adjacency matrix. If }$G$\emph{ contains no
}$C_{2k+2},$\emph{ then}
\[
\mu^{2}-k\mu\leq k\left(  n-1\right)  .
\]
\medskip

\section{Main results}

We write $\left\vert X\right\vert $ for the cardinality of a finite set $X.$
Let $G$ be a graph$,$ and $X$ and $Y$ be disjoint sets of vertices of $G.$ We
write:\medskip

- $V\left(  G\right)  $ for the vertex set of $G$ and $\left\vert G\right\vert
$ for $\left\vert V\left(  G\right)  \right\vert ;$

- $e_{G}\left(  X\right)  $ for the number of edges induced by $X;$

- $e_{G}\left(  X,Y\right)  $ for the number of edges joining vertices in $X$
to vertices in $Y;$

- $G-u$ for the graph obtained by removing the vertex $u\in V\left(  G\right)
;$

- $\Gamma_{G}\left(  u\right)  $ for the set of neighbors of the vertex $u$
and $d_{G}\left(  u\right)  $ for $\left\vert \Gamma_{G}\left(  u\right)
\right\vert .$\medskip

The main result of this note is the following lemma.

\begin{lemma}
\label{le1} Suppose that $k\geq1$ and let the vertices of a graph $G$ be
partitioned into two sets $A$ and $B$.

(A) If
\begin{equation}
2e_{G}\left(  A\right)  +e_{G}\left(  A,B\right)  >\left(  2k-2\right)
\left\vert A\right\vert +k\left\vert B\right\vert , \label{cond1}%
\end{equation}
then there exists a path of order $2k$ or $2k+1$ with both ends in $A.$

(B) If
\begin{equation}
2e_{G}\left(  A\right)  +e_{G}\left(  A,B\right)  >\left(  2k-1\right)
\left\vert A\right\vert +k\left\vert B\right\vert , \label{cond}%
\end{equation}
then there exists a path of order $2k+1$ with both ends in $A.$
\end{lemma}

Note that if we choose the set $B$ to be empty, Lemma \ref{le1} amounts to a
classical result of Erd\H{o}s and Gallai:\medskip

\emph{If a graph of order }$n$\emph{ has more than }$kn/2$\emph{ edges, then
it contains a path of order }$k+2.$\medskip

We postpone the proof of Lemma \ref{le1} to Section \ref{secp} and turn to two consequences.

\begin{theorem}
\label{th1}Let $G$ be a graph with $n$ vertices and $m$ edges. If $G$ does not
contain a$\ C_{2k+2},$ then
\[
\sum_{u\in V\left(  G\right)  }d_{G}^{2}\left(  u\right)  \leq2km+k\left(
n-1\right)  n.
\]

\end{theorem}

\begin{proof}
Let $u$ be any vertex of $G$. Partition the vertices of the graph $G-u$ into
the sets $A=\Gamma_{G}\left(  u\right)  $ and $B=V\left(  G\right)
\backslash\left(  \Gamma_{G}\left(  u\right)  \cup\left\{  u\right\}  \right)
.$ Since $G$ contains no $C_{2k+2},$ the graph $G-u$ dos not contain a path of
order $2k+1$ with both ends in $A.$ Applying Lemma \ref{le1}, part (B), we see
that
\[
2e_{G-u}\left(  A\right)  +e_{G-u}\left(  A,B\right)  \leq\left(  2k-1\right)
\left\vert A\right\vert +k\left\vert B\right\vert ,
\]
and therefore,
\begin{align*}
\sum_{v\in\Gamma_{G}\left(  u\right)  }\left(  d_{G}\left(  v\right)
-1\right)   &  =\sum_{v\in\Gamma_{G}\left(  u\right)  }d_{G-u}\left(
v\right)  =2e_{G-u}\left(  A\right)  +e_{G-u}\left(  A,B\right) \\
&  \leq\left(  2k-1\right)  \left\vert A\right\vert +k\left\vert B\right\vert
\\
&  =\left(  2k-1\right)  d_{G}\left(  u\right)  +k\left(  n-d_{G}\left(
u\right)  -1\right)  .
\end{align*}
Rearranging both sides, we obtain%
\[
\sum_{v\in\Gamma_{G}\left(  u\right)  }d_{G}\left(  v\right)  \leq
kd_{G}\left(  u\right)  +k\left(  n-1\right)  .
\]
Adding these inequalities for all vertices $u\in V\left(  G\right)  ,$ we find
out that
\[
\sum_{u\in V\left(  G\right)  }\sum_{v\in\Gamma_{G}\left(  u\right)  }%
d_{G}\left(  v\right)  \leq k\sum_{u\in V\left(  G\right)  }d_{G}\left(
u\right)  +k\left(  n-1\right)  n=2km+k\left(  n-1\right)  n.
\]
To complete the proof of the theorem note that the term $d_{G}\left(
v\right)  $ appears in the left-hand sum exactly $d_{G}\left(  v\right)  $
times, and so
\[
\sum_{u\in V\left(  G\right)  }\sum_{v\in\Gamma_{G}\left(  u\right)  }%
d_{G}\left(  v\right)  =\sum_{v\in V\left(  G\right)  }d_{G}^{2}\left(
v\right)  .
\]

\end{proof}

Here is a corollary of Theorem \ref{th1} that gives the upper bound for the
proof of (\ref{YC}).

\begin{corollary}
\label{cor1}Let $G$ be a graph with $n$ vertices. If $G$ does not contain
a$\ C_{2k+2},$ then for every $p\geq2,$%
\[
\sum_{u\in V\left(  G\right)  }d_{G}^{p}\left(  u\right)  \leq kn^{p}+O\left(
n^{p-1/2}\right)  .
\]

\end{corollary}

\begin{proof}
Letting $m$ be the number of edges of $G,$ we first deduce an upper bound on
$m.$ Theorem \ref{th1} and the AM-QM inequality imply that
\[
\frac{4m^{2}}{n}\leq\sum_{u\in V\left(  G\right)  }d_{G}^{2}\left(  u\right)
\leq2km+k\left(  n-1\right)  n,
\]
and so,%
\begin{equation}
m\leq-kn+n\sqrt{k\left(  n-1\right)  +k^{2}}<n\sqrt{kn}. \label{ubm}%
\end{equation}
Note that much stronger upper bounds on $m$ are known (e.g., see \cite{BoSi74}
and \cite{Ver00}), but this one is simple and unconditional.

Now Theorem \ref{th1} and inequality (\ref{ubm}) imply that
\begin{align*}
\sum_{u\in V\left(  G\right)  }d_{G}^{p}\left(  u\right)   &  <\sum_{u\in
V\left(  G\right)  }n^{p-2}d_{G}^{2}\left(  u\right)  <kn^{p}+2kmn^{p-2}%
<kn^{p}+2\left(  kn\right)  ^{3/2}n^{p-2}\\
&  =kn^{p}+O\left(  n^{p-1/2}\right)  ,
\end{align*}
completing the proof.
\end{proof}

\section{\label{secp}Proof of Lemma \ref{le1}}

To simplify the proof of Lemma \ref{le1} we state two routine lemmas whose
proofs are given only for the sake of completeness.

\begin{lemma}
\label{le2}Let $P=\left(  v_{1},\ldots,v_{p}\right)  $ be a path of maximum
order in a connected non-Hamiltonian graph $G$. Then $p\geq d_{G}\left(
v_{1}\right)  +d_{G}\left(  v_{p}\right)  +1$.
\end{lemma}

\begin{proof}
Indeed, since $P$ is of maximum order, we sse that $\Gamma_{G}\left(
v_{1}\right)  \subset\left\{  v_{1},\ldots,v_{p}\right\}  $ and $\Gamma
_{G}\left(  v_{p}\right)  \subset\left\{  v_{1},\ldots,v_{p}\right\}  .$ Let
\begin{align*}
r  &  =d_{G}\left(  v_{1}\right)  ,\text{ }s=d_{G}\left(  v_{p}\right)
,\text{ }\\
\Gamma_{G}\left(  v_{1}\right)   &  =\left\{  v_{i_{1}},\ldots v_{i_{r}%
}\right\}  ,\text{ }\Gamma_{G}\left(  v_{p}\right)  =\left\{  v_{j_{1}},\ldots
v_{j_{s}}\right\}  .
\end{align*}
Here we assume that
\[
1<i_{1}<\cdots<i_{r}\leq p,\text{ \ \ \ }1\leq j_{1}<\cdots<j_{s}<p.
\]
If $v_{p}$ is joined to $v_{i_{s}-1}$ for some $1\leq s\leq r,$ then the
sequence
\[
\left(  v_{1},\ldots,v_{i_{s}-1},v_{p},v_{p-1},\ldots,v_{i_{s}},v_{1}\right)
\]
is a cycle of order $p.$ Since $G$ is non-Hamiltonian and connected, there is
an edge joining some of the vertices $v_{1},\ldots,v_{p}$ to a vertex in
$V\left(  G\right)  \backslash\left\{  v_{1},\ldots,v_{p}\right\}  .$ Then we
easily obtain a path longer than $P,$ which contradicts the choice of $P.$

Therefore, $v_{p}$ is not connected to any of the vertices $v_{i_{1}-1}%
,\ldots,v_{i_{r}-1}.$ Thus $\left\{  j_{1},\ldots,j_{s}\right\}  $ and
$\left\{  i_{1}-1,\ldots,i_{r}-1\right\}  $ are disjoint subsets of $\left\{
1,\ldots,p-1\right\}  ,$ implying that
\[
p-1\geq r+s=d_{G}\left(  v_{1}\right)  +d_{G}\left(  v_{p}\right)  ,
\]
and completing the proof.
\end{proof}

\begin{lemma}
\label{le3}Let $P=\left(  v_{1},\ldots,v_{p}\right)  $ be a path of maximum
order in a graph $G$. Then either $v_{1}$ is joined to two consecutive
vertices of $P$ or $G$ contains a cycle of order at least $2d_{G}\left(
v_{1}\right)  .$
\end{lemma}

\begin{proof}
Since $P$ is of maximum order, $\Gamma_{G}\left(  v_{1}\right)  \subset
\left\{  v_{1},\ldots,v_{p}\right\}  .$ Let $\left\{  v_{i_{1}},\ldots
v_{i_{r}}\right\}  =\Gamma_{G}\left(  v_{1}\right)  ,$ where
\[
1<i_{1}<\cdots<i_{r}\leq p.
\]
Assume $v_{1}$ is not joined to two consecutive vertices of $P$, that is to
say, $i_{t}-i_{t-1}\geq2$ for every $t=2,\ldots,r.$ Then the sequence
\[
\left(  v_{1},v_{i_{1}},v_{i_{i}+1},\ldots,v_{i_{r}-1},v_{i_{r}},v_{1}\right)
\]
is a cycle of order at least $1+r+r-1=2r=2d_{G}\left(  v_{1}\right)  ,$
completing the proof.
\end{proof}

\begin{proof}
[\textbf{Proof of Lemma \ref{le1}}]For convenience we shall assume that the
set $B$ is independent. Also, we shall call a path with both ends in $A$ an
$A$\emph{-path}.

\begin{claim}
\label{cl0}If $G$ contains an $A$-path of order $p>2,$ then $G$ contains an
$A$-path of order $p-2.$
\end{claim}

Indeed, let $\left(  v_{1},\ldots,v_{p}\right)  $ be an $A$-path. If $v_{2}\in
B,$ then $v_{3}\in A,$ and so $\left(  v_{3},\ldots,v_{p}\right)  $ is an
$A$-path of order $p-2.$ If $v_{p-1}\in B,$ then $v_{p-2}\in A,$ and so
$\left(  v_{1},\ldots,v_{p-2}\right)  $ is an $A$-path of order $p-2.$
Finally, if both $v_{2}\in A$ and $v_{p-1}\in A,$ then $\left(  v_{2}%
,\ldots,v_{p-1}\right)  $ is an $A$-path of order $p-2.$

The proofs of the two parts of Lemma \ref{le1} are very similar, but since
they differ in the details, we shall present them separately. \medskip

\textbf{Proof of part (A)}\medskip

From Claim \ref{cl0} we easily obtain the following consequence:

\begin{claim}
\label{cl3} If $G$ contains an $A$-path of order $p\geq2k,$ then $G$ contains
an $A$-path of order $2k$ or $2k+1.$
\end{claim}

This in turn implies

\begin{claim}
\label{cl4}If $G$ contains a cycle $C_{p}$ for some $p\geq2k+1,$ then $G$
contains an $A$-path of order $2k$ or $2k+1.$
\end{claim}

Indeed, let $C=\left(  v_{1},\ldots,v_{p},v_{1}\right)  $ be a cycle of order
$p\geq2k+1$. The assertion is obvious if $C$ is entirely in $A,$ so let assume
that $C$ contains a vertex of $B,$ say $v_{1}\in B.$ Then $v_{2}\in A$ and
$v_{p}\in A;$ hence, $\left(  v_{2},\ldots,v_{p}\right)  $ is an $A$-path of
order at least $2k.$ In view of Claim \ref{cl3}, this completes the proof of
Claim \ref{cl4}.

To complete the proof of part (A) we shall use induction on the order of $G.$
First we show that condition (\ref{cond1}) implies that $\left\vert
G\right\vert \geq2k.$ Indeed, assume that $\left\vert G\right\vert \leq2k-1.$
We have
\[
\left\vert A\right\vert ^{2}-\left\vert A\right\vert +\left\vert A\right\vert
\left\vert B\right\vert \geq2e_{G}\left(  A\right)  +e_{G}\left(  A,B\right)
>\left(  2k-2\right)  \left\vert A\right\vert +k\left\vert B\right\vert
\]
and so,
\[
\left\vert G\right\vert \left(  \left\vert A\right\vert -k\right)  =\left(
\left\vert A\right\vert +\left\vert B\right\vert \right)  \left(  \left\vert
A\right\vert -k\right)  >\left(  k-1\right)  \left\vert A\right\vert .
\]
Hence, we find that
\[
\left(  2k-1\right)  \left(  \left\vert A\right\vert -k\right)  >\left(
k-1\right)  \left\vert A\right\vert
\]
and so, $\left\vert A\right\vert >2k-1,$ a contradiction with $\left\vert
A\right\vert \leq\left\vert G\right\vert $.

The conclusion of Lemma \ref{le1}, part (A) follows when $\left\vert
G\right\vert \leq2k-1$ since then the hypothesis is false. Assume now that
$\left\vert G\right\vert \geq2k$ and that the Lemma holds for graphs with
fewer vertices than $G.$ This assumption implies the assertion if $G$ is
disconnected, so to the end of the proof we shall assume that $G$ is connected.

We can assume that $G$ is non-Hamiltonian. Indeed, in view of Claim \ref{cl4},
this is obvious when $\left\vert G\right\vert >2k.$ If $\left\vert
G\right\vert =2k$ and $G$ is Hamiltonian, then no two consecutive vertices
along the Hamiltonian cycle belong to $A,$ and since $B$ is independent, we
have $\left\vert B\right\vert =\left\vert A\right\vert =k.$ Then
\[
k\left(  2k-1\right)  \geq2e_{G}\left(  A\right)  +e_{G}\left(  A,B\right)
>\left(  2k-2\right)  \left\vert A\right\vert +k\left\vert B\right\vert
=k\left(  2k-1\right)  ,
\]
contradicting (\ref{cond1}). Thus, we shall assume that $G$ is non-Hamiltonian.

The induction step is completed if there is a vertex $u\in B$ such that
$d_{G}\left(  u\right)  \leq k.$ Indeed the sets $A$ and $B^{\prime
}=B\backslash\left\{  u\right\}  $ partition the vertices of $G-u$ and also%
\begin{align*}
2e_{G-u}\left(  A\right)  +e_{G-u}\left(  A,B\right)   &  =2e_{G}\left(
A\right)  +e_{G}\left(  A,B\right)  -d_{G}\left(  u\right)  >\left(
2k-2\right)  \left\vert A\right\vert +k\left\vert B\right\vert -k\\
&  =\left(  2k-2\right)  \left\vert A\right\vert +k\left\vert B^{\prime
}\right\vert ;
\end{align*}
hence $G-u$ contains an $A$-path of order at least $2k,$ completing the proof.
Thus, to the end of the proof we shall assume that\medskip\ 

\emph{(a) }$d_{G}\left(  u\right)  \geq k+1$\emph{ for every vertex }$u\in
B.$\medskip

For every vertex $u\in A,$ write $d_{G}^{\prime}\left(  u\right)  $ for its
neighbors in $A$ and $d_{G}^{\prime\prime}\left(  u\right)  $ for its
neighbors in $B.$ The induction step can be completed if there is a vertex
$u\in A$ such that $2d_{G}^{\prime}\left(  u\right)  +d_{G}^{\prime\prime
}\left(  u\right)  \leq2k-2.$ Indeed, if $u$ is such a vertex, note that the
sets $A^{\prime}=A\backslash\left\{  u\right\}  $ and $B$ partition the
vertices of $G-u$ and also%
\begin{align*}
2e_{G-u}\left(  A\right)  +e_{G-u}\left(  A,B\right)   &  =2e_{G}\left(
A\right)  +e_{G}\left(  A,B\right)  -2d_{G}^{\prime}\left(  u\right)
-d_{G}^{\prime\prime}\left(  u\right) \\
&  >\left(  2k-2\right)  \left\vert A\right\vert +k\left\vert B\right\vert
-2k+2\\
&  =\left(  2k-2\right)  \left\vert A^{\prime}\right\vert +k\left\vert
B\right\vert ;
\end{align*}
hence $G-u$ contains an $A$-path of order at least $2k,$ completing the proof.
Hence we have $2d_{G}^{\prime}\left(  u\right)  +d_{G}^{\prime\prime}\left(
u\right)  \geq2k-1,$ and so $d_{G}\left(  u\right)  \geq k.$ Thus, to the end
of the proof, we shall assume that:\medskip

\emph{(b) }$d_{G}\left(  u\right)  \geq k$ \emph{for every vertex }$u\in
A.$\medskip

Select now a path $P=\left(  v_{1},\ldots,v_{p}\right)  $ of maximum length in
$G.$ To complete the induction step we shall consider three cases: \emph{(i)
}$v_{1}\in B,$ $v_{p}\in B;$ \emph{(ii)} $v_{1}\in B,$ $v_{p}\in A,$ and
\emph{(iii) }$v_{1}\in A,$ $v_{p}\in A.$\medskip

\emph{Case (i): }$v_{1}\in B,$ $v_{p}\in B$\medskip

In view of assumption \emph{(a)} we have $d_{G}\left(  v_{1}\right)
+d_{G}\left(  v_{p}\right)  \geq2k+2,$ and Lemma \ref{le2} implies that
$p\geq2k+3$. We see that $\left(  v_{2},\ldots,v_{p-1}\right)  $ is an
$A$-path of order at least $2k+1,$ completing the proof by Claim \ref{cl3}.
\medskip

\emph{Case (ii): }$v_{1}\in B,$ $v_{p}\in A$\medskip

In view of assumptions \emph{(a) }and \emph{(b)} we have $d_{G}\left(
v_{1}\right)  +d_{G}\left(  v_{p}\right)  \geq2k+1,$ and Lemma \ref{le2}
implies that $p\geq2k+2,$ and so, $\left(  v_{2},\ldots,v_{p}\right)  $ is an
$A$-path of order at least $2k+1$. This completes the proof by Claim
\ref{cl3}.\medskip

\emph{Case (iii): }$v_{1}\in A,$ $v_{p}\in A$\medskip

In view of assumption \emph{(b) }we have $d_{G}\left(  v_{1}\right)
+d_{G}\left(  v_{p}\right)  \geq2k,$ and Lemma \ref{le2} implies that
$p\geq2k+1.$ Since $\left(  v_{1},\ldots,v_{p}\right)  $ is an $A$-path of
order at least $2k+1,$ by Claim \ref{cl3}, the proof of part (A) of Lemma
\ref{le1} is completed.\medskip

\textbf{Proof of part (B)}\medskip

From Claim \ref{cl0} we easily obtain the following consequence:

\begin{claim}
\label{cl1} If $G$ contains an $A$-path of odd order $p\geq2k+1,$ then $G$
contains an $A$-path of order exactly $2k+1.$
\end{claim}

From Claim \ref{cl1} we deduce another consequence:

\begin{claim}
\label{cl2}If $G$ contains a cycle $C_{p}$ for some $p\geq2k+1,$ then $G$
contains an $A$-path of order exactly $2k+1.$
\end{claim}

Indeed, let $C=\left(  v_{1},\ldots,v_{p},v_{1}\right)  $ be a cycle of order
$p\geq2k+1$. If $p$ is odd, then some two consecutive vertices of $C$ belong
to $A,$ say the vertices $v_{1}$ and $v_{2}.$ Then $\left(  v_{2},\ldots
,v_{p},v_{1}\right)  $ is an $A$-path of odd order $p\geq2k+1,$ and by Claim
\ref{cl1} the assertion follows. If $p$ is even, then $p\geq2k+2$. The
assertion is obvious if $C$ is entirely in $A,$ so let assume that $C$
contains a vertex of $B,$ say $v_{1}\in B.$ Then $v_{2}\in A$ and $v_{p}\in
A;$ hence $\left(  v_{2},\ldots,v_{p}\right)  $ is an $A$-path of odd order at
least $2k+1,$ completing the proof of Claim \ref{cl2}.

To complete the proof of Lemma \ref{le1} we shall use induction on the order
of $G.$ First we show that condition (\ref{cond}) implies that $\left\vert
G\right\vert \geq2k+1.$ Indeed, assume that $\left\vert G\right\vert \leq2k.$
We have
\[
\left\vert A\right\vert ^{2}-\left\vert A\right\vert +\left\vert A\right\vert
\left\vert B\right\vert \geq2e_{G}\left(  A\right)  +e_{G}\left(  A,B\right)
>\left(  2k-1\right)  \left\vert A\right\vert +k\left\vert B\right\vert
\]
and so,
\[
\left\vert G\right\vert \left(  \left\vert A\right\vert -k\right)  =\left(
\left\vert A\right\vert +\left\vert B\right\vert \right)  \left(  \left\vert
A\right\vert -k\right)  >k\left\vert A\right\vert .
\]
Hence, we find that $2k\left(  \left\vert A\right\vert -k\right)  >k\left\vert
A\right\vert ,$ and $\left\vert A\right\vert >2k,$ contradicting that
$\left\vert A\right\vert \leq\left\vert G\right\vert .$

The conclusion of Lemma \ref{le1}, part (B) follows when $\left\vert
G\right\vert \leq2k$ since then the hypothesis is false. Assume now that
$\left\vert G\right\vert \geq2k+1$ and that the assertion holds for graphs
with fewer vertices than $G.$ This assumption implies the assertion if $G$ is
disconnected, so to the end of the proof we shall assume that $G$ is
connected. Also, in view of Claim \ref{cl2} and $\left\vert G\right\vert
\geq2k+1,$ we shall assume that $G$ is non-Hamiltonian.

The induction step is completed if there is a vertex $u\in B$ such that
$d_{G}\left(  u\right)  \leq k.$ Indeed the sets $A$ and $B^{\prime
}=B\backslash\left\{  u\right\}  $ partition the vertices of $G-u$ and also%
\begin{align*}
2e_{G-u}\left(  A\right)  +e_{G-u}\left(  A,B\right)   &  =2e_{G}\left(
A\right)  +e_{G}\left(  A,B\right)  -d_{G}\left(  u\right) \\
&  >\left(  2k-1\right)  \left\vert A\right\vert +k\left\vert B\right\vert
-k\\
&  =\left(  2k-1\right)  \left\vert A\right\vert +k\left\vert B^{\prime
}\right\vert ;
\end{align*}
hence $G-u$ contains an $A$-path of order $2k+1,$ completing the proof. Thus,
to the end of the proof we shall assume that: \medskip

\emph{(a)} $d_{G}\left(  u\right)  \geq k+1$\emph{ for every vertex }$u\in
B.$\medskip

For every vertex $u\in A,$ write $d_{G}^{\prime}\left(  u\right)  $ for its
neighbors in $A$ and $d_{G}^{\prime\prime}\left(  u\right)  $ for its
neighbors in $B.$ The induction step can be completed if there is a vertex
$u\in A$ such that $2d_{G}^{\prime}\left(  u\right)  +d_{G}^{\prime\prime
}\left(  u\right)  \leq2k-1.$ Indeed, if $u$ is such a vertex, note that the
sets $A^{\prime}=A\backslash\left\{  u\right\}  $ and $B$ partition the
vertices of $G-u$ and also%
\begin{align*}
2e_{G-u}\left(  A\right)  +e_{G-u}\left(  A,B\right)   &  =2e_{G}\left(
A\right)  +e_{G}\left(  A,B\right)  -2d_{G}^{\prime}\left(  u\right)
-d_{G}^{\prime\prime}\left(  u\right) \\
&  >\left(  2k-1\right)  \left\vert A\right\vert +k\left\vert B\right\vert
-2k+1\\
&  =\left(  2k-1\right)  \left\vert A^{\prime}\right\vert +k\left\vert
B\right\vert ;
\end{align*}
hence $G-u$ contains an $A$-path of order $2k+1,$ completing the proof. Thus,
to the end of the proof, we shall assume that:\medskip

\emph{(b) }$d_{G}\left(  u\right)  \geq k$ \emph{for every vertex }$u\in A$
\emph{and if }$u$\emph{ has neighbors in }$B,$\emph{ then }$d_{G}\left(
u\right)  \geq k+1.$\medskip

Select now a path $P=\left(  v_{1},\ldots,v_{p}\right)  $ of maximum length in
$G.$ To complete the induction step we shall consider three cases: \emph{(i)
}$v_{1}\in B,$ $v_{p}\in B;$ \emph{(ii)} $v_{1}\in B,$ $v_{p}\in A,$ and
\emph{(iii) }$v_{1}\in A,$ $v_{p}\in A.$\medskip

\emph{Case (i): }$v_{1}\in B,$ $v_{p}\in B$\medskip

In view of assumption \emph{(b) }we have $d_{G}\left(  v_{1}\right)
+d_{G}\left(  v_{p}\right)  \geq2k+2,$ and Lemma \ref{le2} implies that
$p\geq2k+3$. If $p$ is odd, we see that $\left(  v_{2},\ldots,v_{p-1}\right)
$ is an $A$-path of order at least $2k+1,$ and by Claim \ref{cl1}, the proof
is completed.

Suppose now that $p$ is even. Applying Lemma \ref{le3}, we see that either $G$
has a cycle of order at least $2d_{G}\left(  v_{1}\right)  \geq2k+2,$ or
$v_{1}$ is joined to $v_{i}$ and $v_{i+1}$ for some $i\in\left\{
2,\ldots,p-2\right\}  .$ In the first case we complete the proof by Claim
\ref{cl2}; in the second case we see that the sequence
\[
\left(  v_{2},v_{3},\ldots,v_{i},v_{1},v_{i+1},v_{i+2},\ldots,v_{p-1}\right)
\]
is an $A$-path of order $p-1.$ Since $p-1$ is odd and $p-1\geq2k+3,$ the proof
is completed by Claim \ref{cl1}.\medskip

\emph{Case (ii): }$v_{1}\in B,$ $v_{p}\in A$\medskip

In view of assumptions \emph{(a)} and \emph{(b)} we have $d_{G}\left(
v_{1}\right)  +d_{G}\left(  v_{p}\right)  \geq2k+1,$ and Lemma \ref{le2}
implies that $p\geq2k+2$. If $p$ is even, we see that $\left(  v_{2}%
,\ldots,v_{p-1}\right)  $ is an $A$-path of order at least $2k+1,$ and by
Claim \ref{cl1}, the proof is completed.

Suppose now that $p$ is odd. Applying Lemma \ref{le3}, we see that either $G$
has a cycle of order at least $2d_{G}\left(  v_{1}\right)  \geq2k+2,$ or
$v_{1}$ is joined to $v_{i}$ and $v_{i+1}$ for some $i\in\left\{
2,\ldots,p-1\right\}  .$ In the first case we complete the proof by Claim
\ref{cl2}; in the second case we see that the sequence
\[
\left(  v_{2},v_{3},\ldots,v_{i},v_{1},v_{i+1},v_{i+2},\ldots,v_{p}\right)
\]
is an $A$-path of order $p.$ Since $p$ is odd and $p\geq2k+2,$ the proof is
completed by Claim \ref{cl1}.\medskip

\emph{Case (iii): }$v_{1}\in A,$ $v_{p}\in A$\medskip

In view of assumption \emph{(b)} we have $d_{G}\left(  v_{1}\right)
+d_{G}\left(  v_{p}\right)  \geq2k,$ and Lemma \ref{le2} implies that
$p\geq2k+1.$ If $p$ is odd, the proof is completed\ by Claim \ref{cl1}.

Suppose now that $p$ is even, and therefore, $p\geq2k+2.$ If $v_{2}\in A,$
then the sequence $\left(  v_{2},\ldots,v_{p}\right)  $ is an $A$-path of odd
order $p-1\geq2k+1,$ completing the proof by Claim \ref{cl1}. If $v_{2}\in B,$
we see that $v_{1}$ has a neighbor in $B,$ and so, $d_{G}\left(  v_{1}\right)
\geq k+1.$

Applying Lemma \ref{le3}, we see that either $G$ has a cycle of order at least
$2d_{G}\left(  v_{1}\right)  \geq2k+2,$ or $v_{1}$ is joined to $v_{i}$ and
$v_{i+1}$ for some $i\in\left\{  2,\ldots,p-2\right\}  .$ In the first case we
complete the proof by Claim \ref{cl2}. In the second case we shall exhibit an
$A$-path of order $p-1.$ Indeed, if $i=2,$ let
\[
Q=\left(  v_{1},v_{3},v_{4},\ldots,v_{p}\right)  ,
\]
and if $i\geq3,$ let
\[
Q=\left(  v_{3},\ldots,v_{i},v_{1},v_{i+1},v_{i+2},\ldots,v_{p}\right)  .
\]
In either case $Q$ is an $A$-path of order $p-1.$ Since $p-1$ is odd and
$p-1\geq2k+1,$ the proof is completed by Claim \ref{cl1}.

This completes the proof of Lemma \ref{le1}.
\end{proof}

\textbf{Acknowledgment }Thanks are due to Dick Schelp and Ago Riet for useful
discussions on Lemma \ref{le1}.


\bigskip

\begin{thebibliography}{9}                                                                                                %


\bibitem {Bol98}B. Bollob\'{a}s, \emph{Modern Graph Theory,} Graduate Texts in
Mathematics, \textbf{184,} Springer-Verlag, New York (1998), xiv+394 pp.

\bibitem {BoSi74}J. A. Bondy and M. Simonovits, Cycles of even length in
graphs, \emph{J. Comb. Theory Ser. B} \textbf{16} (1974), 97--105.

\bibitem {CaYu00}Y. Caro, R. Yuster, A Tur\'{a}n type problem concerning the
powers of the degrees of a graph, \emph{Electron. J. Comb.} \textbf{7 }(2000),
RP 47.

\bibitem {ErGa59}P. Erd\H{o}s, T. Gallai, On maximal paths and circuits of
graphs, \emph{Acta Math. Acad. Sci. Hungar.} \textbf{10} (1959), 337--356.

\bibitem {Nik08}V. Nikiforov, The spectral radius of graphs with forbidden
paths and cycles, \emph{preprint.}

\bibitem {Ver00}J. Verstra\"{e}te, On arithmetic progressions of cycle lengths
in graphs, \emph{Combin. Probab. Comput.} \textbf{9} (2000), 369--373.
\end{thebibliography}
\end{document}